\magnification=1200
\input amssym.def
\input amssym.tex
 \font\newrm =cmr10 at 24pt
\def\bul{\raise .9pt\hbox{\newrm .\kern-.105em } }

 \def\fr{\frak}

 \baselineskip 20pt
 
 \def\h{\hbox{ }}

 \def\r{{\fr r}}

 \def\n{{\fr n}}

 \def\ss{{\fr s}}
 \def\k{{\fr k}}
 \def\b{{\fr b}}
 
 \def\hh{{\fr h}}

 \def\g{{\fr g}}

 \def\<{\le}
 \def\>{\ge}

 \def\s{{\h\subset\h}}
 
 \def\vs{\vskip }

 \def\mapright#1
  {\smash{\mathop
  {\longrightarrow}
  \limits^{#1}}}

 \def\kk#1{{\kern .4 em} #1}
 \def\vs{\vskip 1pc}

\centerline{\bf Generalized Amitsur--Levitski theorem and equations for Sheets}
 \centerline{\bf in  a Reductive  complex Lie algebra} \vskip1.5pc
 
 \centerline{ \it Bertram Kostant}\vskip 1.5pc
 
 \centerline{\bf  0. Introduction}\vskip 1.5pc 
 
 {\bf 0.1.}The Amitsur-Levitski theorem is a famous result stating that matrix algebras satisfy a certain identity often referred to as the
  standard identity. This result can be generalized when viewed as a statement in Lie theory. If $\g$ is a complex reductive Lie algebra the generalization involves the nilcone of $\g$. 
 
 If $U$ is the universal enveloping algebra then  in [2] we have shown  $$U= Z\otimes E\eqno(0.1)$$ where $E$  is the graded $\g$--module (under the adjoint action) spanned all powers of nilpotent elements in $\g$ and $Z=Cent\,U$. Since $U$ is a quotient of the tensor algebra of $\g$ one has a $\g$-map $$\Gamma_T:A(\g)\to U\eqno(0.2)$$ where $A(\g)$ is the graded  $\g$-module of alternating tensors. We have proved in (4) that $$\Gamma_T(A^{2k}(\g))\s E^{k}\eqno (0.3)$$ for any $k\in \Bbb Z$. If $\pi$ is a representation of $U$ on a vector space $V$ this obviously implies $$Ker\, \pi|\Gamma_T (A^{2k}(\g))\s Ker\,\pi|E^{k}\eqno (0.4)$$ This as one easily notes generalizes the Amitsur--Levitski theorem. \vs {0.2.} Since the exterior algebra $\wedge \g$ is also a quotient of the tensor algebra the quotient map defines an isomorphism $$A(\g)\to \wedge \g$$ of graded $\g$-modules.  In particular by restriction an isomorphism  $$A^{even}(\g) \to \wedge^{even}\g\eqno (0.5)$$ noting the right side of (0.5) is a commutative algebra. 
 
 Identify $\g$ with its dual $\g^*$ using the Killing form. We now consider the ``commutative" analogue of the statements in \S 0.1. Let $P(\g)$ be the ring of polynomial functions on $\g$. Analogous to (0.1) one has $$P(\g) = J\otimes H\eqno (0.6)$$ where $J$ is the algebra $P(\g)^{\g}$ of polynomial invariants on $\g$ and $H$ is the graded $\g$-module of harmonic polynomials on $g$. (See Theorem 11 in [2]). Poincar\'e--Birkhoff--Witt symmetrization induces a $\g$-module isomorphism $$U\to P(\g)\eqno (0.7)$$ which restricts to a $\g$-module isomorphisms $$\eqalign{Z& \to J \cr E&\to H\cr}\eqno (0.8)$$
 
 Using these isomorphisms (0.2) and (0.3) now define $\g$-module maps $$\Gamma: \wedge^{even}\g \to P(\g)\eqno (0.9)$$ and for $\k\in \Bbb Z$,  $$\Gamma:\wedge^{2k}\g\to H^k\eqno (0.10 )$$\vskip .5pc Let $R^k(\g)$ be the image of (0.10) so that $R^k(\g)$ is a $\g$-module of harmonic polynomials on $\g$ of degree $k$. The significance of $R^k(\g)$ has to do with the dimensions of $Ad\,\g$ adjoint (=coadjoint) orbits. Any such orbit is symplectic and hence is even dimensional. For $j\in \Bbb Z$ let $\g^{(2j)}= \{x\in \g\mid dim\,[\g,x] = 2j\}$. We recall that a $2j\,\g$-sheet is an irreducible component of $\g^{(2j)}$. Let $Var\,R^k(\g) = \{x\in \g\mid p(x)=0,\,\forall p\in R^k(\g)\}$. In [4] (see Proposition 3.2 in [4])  we prove \vs {\bf Theorem 0.1. } {\it One has $$Var\,R^k(\g) = \cup_{2j<2k}\, \g^{(2j)}\eqno (0.11)$$ or that $Var\,R^k(\g)$ is the set of all $2j\,\g$-sheets for $j<k$.} \vs One can explicitly describe the $\g$-module $R^k(\g)$. (See \S 1.2 in [4]). 
 
 One notes that that there exists a subset $\Pi(k)$ of the symmetric group $sym(k)$, having  the cardinality of is $(2k-1)(2k-3)\cdots 1$, such that  the correspondence $$\nu \mapsto ((\nu(1),\nu(2)), (\nu(3),\nu(4)),\ldots ,(\nu((2k-1)),\nu(2k))) $$ sets up a bijection of $\Pi(k)$ with the set of all partitions of $(1,2,\ldots,2k)$ into a union of subsets each  of which has two elements. We also observe that $\Pi(k)$ can and will be chosen such that $sg\,\nu =1$ for all $\nu\in \Pi(k)$. 
 
 The following is a restatement of the results in \S3.2 of [4] (see especially (3.25) and (3.29) in [4])\vs {\bf Theorem 0.2.} {\it For any $k\in \Bbb Z$ there exists a nonzero scalar $c_k$ such that for any  $x_i\,i=1,\dots, 2k,$ in $\g$
 $$\Gamma (x_1\wedge \cdots \wedge x_{2k}) = c_k \sum_{\nu\in \Pi(k)} [x_{\nu(1)},x_{\nu(2)}]\cdots [x_{\nu(2k-1)},x_{\nu(2k)}] \eqno (0.12))$$
 Furthermore the homogeneous polynomial of degree $k$ on the right side of (0.12) is harmonic and $R^k(\g)$ is the span of all such polynomials for an arbitrary choice of the $x_i$.} \vskip 1.5pc \centerline{\bf 0. Introduction to results of our joint paper, [5], with Nolan Wallach}\vskip 1.5 pc {0.2.} Henceforth assume that $\g$ is simple so that  the adjoint representation is irreducible. Let $n=dim\,\g$ and let $\hh$ be a Cartan subalgebra of $\g$ so that $\ell=dim\,\hh$ is the rank of $\g$. Let $\Delta_+$ be a choice of positive roots for $(\g,\hh)$ so that $n= \ell +2r$ where $r= card\,\Delta_+$.  One readily has that $R^k(\g)= 0$ if $k>r$  and $k=r$ is the maximal value of $k$ for which $R^k(\g)$ is interesting. In fact $R^r(\g)$ is the variety of singular elements of $\g$. The paper [5] is devoted to a study of $R^r(\g)$ and in particular to its remarkable $\g$-module structure. Of course Theorem 0.2 may be applied to determine $R^r(\g)$. However in [5]  we establish a different determination of $R^r(\g)$. The results of [2] imply that the adjoint representation, has multiplicity $\ell$ in $H$. Choosing a basis of $\g$ and suitable generators of $J$ this means there exists a $\ell\times n$ matrix $Q$ whose rows are the $\ell$ occurrences of the adjoint representation in $H$. One notes that this matrix has $n\choose \ell$ minors of size $\ell \times \ell$. We prove \vs {\bf Theorem 0.3.} {\it The determinant of any $\ell\times\ell$ minor of $Q$ is an element of $R^r(\g)$ and indeed $R^r(\g) $ is the span of the determinants of all these minors. } \vs If $\varphi \in \Delta_+$ let $e_{\varphi}\in \g$ be a corresponding root vector.  Let $\n$ be the Lie algebra spanned by $e_{\varphi}$ for $\varphi\in \Delta_+$ and let $\b$ the Borel subalgebra of $\g$ defined by putting $\b= \hh + \n$. Now a subset $\Phi\s \Delta_+$ will be called an ideal in $\Delta_+$ if the span, $\n_{\Phi}$, of $e_{\varphi}$, for $\varphi\in \Phi$, is an ideal of $\b$. In such case let  $$\langle \Phi\rangle = \sum_{i=1}^{k} \varphi_i\eqno (0.13) $$ It follows easily that $\langle \Phi\rangle$ is a dominant weight. Let $V_{\Phi}$ be an irreducible $\g$-module having $\langle \Phi\rangle$ as highest weight. In [4] we have shown if $\Phi, \Psi$ are distinct ideals of $\Delta_+ $ then 
 $$V_{\Phi}\, \,\hbox{and}\,V_{\Psi},\,\, \hbox{define inequivalent representations of $\g$}\eqno (0.14)$$\vskip .5pc 
  Let  ${\cal I}(\ell)$ be the set all ideals $\Phi$ in $\Delta_+$ such that $dim\,\n_{\Phi} = \ell$. \vs {\bf Remark 0.1} One notes that if $\g$ is of type $A_{\ell}$ then $card\, {\cal I}(\ell)$ is $P(\ell)$ where $P$ here is the classical partition function 
 \vs Let $Cas\in Z$ be the quadratic Casimir element normalized by the condition that the eigenvalue of $Cas$ in the adjoint representation is 1. The following theorem, giving the rather striking $\g$--module structure of $R^r(\g)$ is one of the main result of [5]. The theory of abelian ideals of $\b$ is an area of considerable current research activity. The proof of theorem makes direct contact with this theory.\vs {\bf Theorem 0.4.} {\it For any $\Phi\in {\cal I}(\ell)$ the $\b$--ideal $\n_{\Phi}$ is abelian. Furthermore as $\g$--modules one has the equivalence  $$R^r(\g) \cong \sum_{\Phi\in {\cal I}(\ell)} V_{\Phi}\eqno (0.15)$$ so that $R^r(\g)$ is a multiplicity one $\g$--module with $card \,{\cal I}(\ell)$ irreducible components. In addition $Cas$ has the eigenvalue $\ell$ on each and every such component}\vskip 1.5pc \centerline{\bf  Results}\vskip 1.5pc

  {\bf 1.} Let $R$ be an associative ring and for any $k\in \Bbb Z$ and $x_i,\ldots, x_k,$ in $R$ one defines an alternating sum of products 
 
 $$[[x_1,\ldots x_k]] = \sum_{\sigma \in Sym\,k} sg(\sigma)\, x_{\sigma(1)} \cdots x_{\sigma(k)}\eqno (1)$$
 
One  says that $R$ satisfies the standard identity of degree $k$ if $[[x_1,\ldots, x_k]] =0$ for any choice of the $ x_i \in R$. Of course $R$ is commutative if  and only if it satisfies the standard identity of degree 2. 

Now for any $n\in \Bbb Z$ and field $F$, let $M(n,F)$ be the algebra of $n \times n$ matrices over $F$. The following is the famous Amitsur--Levitski theorem. \vs

{\bf Theorem 1}. {\it $M(n,F)$ satisfies the standard identity of degree 2n}. \vs {\bf Remark 1}.  By restricting to matrix units , for a proof, it suffices to take $F=\Bbb C$. \vs
Without any knowledge that it was a known theorem  we came upon Theorem 1 in [1],  a long time ago, from the point of Lie algebra cohomology. In fact the result  follows from the fact that if $\g = M(n,\Bbb C)$, then the restriction to $\g$ of the primitive cohomology class of  degree $2n+1$ of $M(n+1,\Bbb C)$ to $\g$ vanishes. 

Of course  $\g_1\s \g$ where $\g_1= Lie \,SO(n,\Bbb C)$. Assume $n$ is even. One proves that the restriction to $\g_1$ of the primitive class of degree $2n-1$ (highest primitive class) of $\g$ vanishes on $\g_1$. This leads to a new standard identity,  namely \vs {\bf Theorem 2}. {\it  $$[[x_1,\dots, x_{2n-2}]] = 0 \eqno (2)$$ for any choice of $x_i\in\g_1$. That is any choice of skew-symmetric matrices. }\vs {\bf Remark 2}. Theorem 2 is immediately evident when  $n=2$.

Theorems 1 and 2 suggest that standard identities can be viewed as a subject in Lie theory. Theorem 3 below  offers support for this idea. Let $\r$ be a complex reductive Lie algebra and let $$\pi:\r \to End \,V\eqno (3) $$ be a finite--dimensional complex  completely reducible representation. If $w\in \r$ is nilpotent then $\pi(w)^k=0$ for some $k\in \Bbb Z$. Let $\varepsilon (\pi)$ be the minimal integer $k$ such that $\pi(w)^k= 0$ for all nilpotent $w\in\r$. In case $\pi$ is irreducible one can easily give a formula for $\varepsilon (\pi)$ in terms of the highest weight. If $\g$ (resp. $\g_1$) is given as above and $\pi$ (resp.  $\pi_1)$ ) is the defining representation 
 then $\varepsilon(\pi) = n$ and $\varepsilon (\pi_1) = n-1$. Consequently the following theorem (see [4])  generalizes Theorems 1 and 2.\vs {\bf Theorem 3.} {\it Let $\r$ be a complex reductive Lie algebra and let $\pi$ be as above. Then for any $x_i\in \r,\,i=1,\ldots, 2\varepsilon(\pi),$ one has  $$[[\hat {x}_1,\ldots, \hat {x}_{2\,\varepsilon(\pi)}]] = 0\eqno (4)$$ where $\hat {x}_i= \pi (x_i)$.} \vs   {\bf 2.} Henceforth $\g$, until mentioned otherwise,  will be an arbitrary reductive complex finite dimensional Lie algebra. Let $T(\g)$ be the tensor  algebra over $\g$ and let $S(\g)\s T(\g)$ ( resp. $A(g)\s T(\g)$ )  be the subspace of symmetric (resp. alternating) tensors in 
 $T(\g)$. The natural grading on $T(\g)$ restricts to a grading on $S(\g)$ and $A(\g)$. In particular, where multiplication is tensor product one notes \vs {\bf Proposition 1.} {\it $A^j(\g)$ is the span of $[[x_1,\ldots ,x_j]]$ over all choices of $x_i,\,i=1,\ldots,j, $ in $\g$. } \vs Now let $U(\g)$ be the universal enveloping algebra of $\g$. Then $U(\g)$ is the quotient algebra of $T(\g)$ so that there is an algebra epimorphism $$\tau :T(\g)\to U(\g)$$ Let $Z = Cent\,U(\g)$ and let $E\s U(g)$ be the graded subspace spanned by all powers $e^j, j=1,...$ where $e\in \g$ is nilpotent. In [2] (see Theorem 21 in [2]) we proved (where tensor product identifies with multiplication) $$U(\g)=
 Z\otimes E\eqno(5)$$ In[4] (see Theorem 3.4 in [4]) we proved \vs {\bf Theorem 4.} {\it  For any $k\in \Bbb Z)$ one has $$\tau(A^{2k}(\g))\s E^k\eqno (6)$$}\vs Theorem 3 is then an immediate consequence of Theorem 4. Indeed, using the notation of Theorem 3,  let $\pi_U: U(\g)\to End \,V$ be the algebra extension of $\pi$ to $U(\g)$. One then has \vs {\bf Theorem 5.} {\it If $E^k\s Ker\,\pi_{U}$ then $$[[\hat {x}_1,\ldots, \hat {x}_{2k}]] = 0\eqno (7)$$ for any $x_i,\ldots, x_{2k}$ in $\g$}. \vs
 
 {\bf 3.} The Poincar\'e--Birkhoff--Witt theorem says that the restriction $\tau:S(\g)\to U(\g)$ is a linear isomorphism.  Consequently, given any $t\in T(\g)$ there exists a unique element $\overline {t}$ in $S(\g)$ such that $$\tau(t) = \tau(\overline{t})\eqno (8)$$ Let $A^{even}(\g)$  be the span of alternating tensors of even degree. Restricting to $A^{even}(\g)$ one has a $\g$-module map $$\Gamma_T:A^{even}(\g) \to S(\g)$$ defined so that if 
 $a \in A^{even}(\g)$ then $$\tau(a) = \tau  (\Gamma_T(a))\eqno (9)$$ Now the (commutative) symmetric algebra $P(\g)$ over $g$ and exterior algebra $\wedge \g$ are quotient algebras of $T(\g)$. The restriction of the quotient map clearly induces $\g$-module  isomorphisms $$\eqalign{\tau_S:S(\g)&\to P(\g)\cr \tau_A:A^{even}(\g)&\to \wedge^{even} \g\cr}\eqno (10)$$ where $ \wedge^{even} \g$ is the commutative subalgebra of $\wedge \g$ spanned by elements of even degree. We may complete the commutative diagram defining $$\Gamma:\wedge^{even} \g \to P(\g) \eqno(11)$$ so that on $A^{even}(\g)$ one has $$\tau_S\circ\Gamma_T= \Gamma \circ \tau_A\eqno (12)$$ By (6) one notes that for $k\in \Bbb Z$ one has $$\Gamma:\wedge^{2k}\g \to P^k(\g) \eqno (13)$$ \vskip .5pc The Killing form extends to a nonsingular symmetric bilinear form on $P(\g)$ and $\wedge \g$. This enables us to identify $P(\g)$ with the algebra of polynomial functions on $\g$ and to identify $\wedge\g$ with its dual space  $\wedge\g^*$ where $\g^*$ is the dual space to $\g$. Let $R^k(\g)$ be the image (13) so that $R^k(\g)$ is a $\g$-module of homogeneous  polynomial functions of degree $k$ on $\g$. 
 The significance of $R^k(\g)$ has to do with the dimensions of $Ad\,\g$ adjoint (=coadjoint) orbits. Any such orbit is symplectic and hence is even dimensional. For $j\in \Bbb Z$ let $\g^{(2j)}= \{x\in \g\mid dim\,[\g,x] = 2j\}$. We recall that a $2j\,\g$-sheet is an irreducible component of $\g^{(2j)}$. Let $Var\,R^k(\g) = \{x\in \g\mid p(x)=0,\,\forall p\in R^k(\g)\}$. In [4] (see Proposition 3.2 in [4])  we prove \vs {\bf Theorem 6. } {\it One has $$Var\,R^k(\g) = \cup_{2j<2k}\, \g^{(2j)}\eqno (14)$$ or that $Var\,R^k(\g)$ is the set of all $2j\,\g$-sheets for $j<k$.}\vs Let $\gamma$ be the transpose of $\Gamma$ . Thus $$\gamma:P(\g)\to \wedge^{even}\g\eqno (15)$$ and one has for $p\in P(\g)$ and $u\in \wedge\g$, $$(\gamma(p),u) = (p,\Gamma(u))\eqno (16)$$ One also notes $$\gamma:P^k(\g) \to \wedge^{2k}\g\eqno (17)$$ A proof of Theorem 6 depends upon establishing some nice algebraic properties of $\gamma$. Since we have, via the Killing form,  identified $\g$ with its dual, $\wedge \g$ is the underlying space for a standard cochain complex $(\wedge\g,d)$  where $d$ is the coboundary operator of degree $+1$.   In particular if $x\in \g$ then $dx \in \wedge^2\g$. Identifying $\g$ here with $P^1(\g)$ one has a map  $$P^1(\g) \to \wedge^2\g \eqno (18)$$ \vs {\bf Theorem 7.} {\it The map (15) is the homomorphism of commutative algebras extending (18). In particular for any $x\in \g$ $$\gamma(x^k) = (-dx)^k\eqno (19)$$}\vs  The  connection with Theorem 6 follows from  \vs {\bf Proposition 2.} {\it Let $x\in \g$. Then $x\in \g^{(2k)}$ if and only if $k$ is maximal such that $(dx)^k\neq 0$, in which case there is a scalar $c\in \Bbb C^{\times}$ such that $$(dx)^k = c\,\, w_1\wedge\cdots \wedge w_{2k}\eqno (20)$$ where $w_i,\,i=1,\ldots,2k$, is a basis of $[x,\g]$. }\vs
 For a proof of Theorem 7 and Proposition 2 see Theorem 1.4 and Proposition 1.3 in [4].
 
  We wish to explicitly describe the $\g$-module $R^k(\g)$. (See \S 1.2 in [4]). Let $J = P(\g)^\g$ so that $J$ is the ring of $ Ad\,\g$ polynomial invariants. Let $Diff\,P(\g)$ be the algebra of differential operators on $P(\g)$ with constant coefficients. One then has an algebra isomorphism $$P(\g)\to Diff\,P(\g), \,\,\, q \mapsto \partial_q$$ where for $p,q,f\in P(\g)$ one has $$(\partial_q p ,f) = (p,qf)\eqno (21)$$ and $\partial_x$, for $x\in \g$, is the partial derivative defined by $x$.

 Let $J_+\s J $ be the $J$-ideal of all $p\in J$ with zero constant term and let $$H= \{q\in P(\g)\mid \partial_pq = 0\,\,\forall p\in J_+\}$$ 
 $H$ is a graded $\g$-module whose elements are called harmonic polynomials. Then one knows (see Theorem 11 in [2]) that, where tensor product is realized by polynomial multiplication, $$P(\g) = J\otimes H\eqno (22)$$ It is immediate from (21) that $H$ is the orthocomplement of the ideal $J_+P(\g)$  in $P(\g)$. However since $\gamma$ is an algebra  homomorphism one has $$J_+P(\g) \s Ker\, \gamma\eqno (23)$$ since one easily has that $J_+\s Ker \,\gamma$. Indeed this is clear since $$\eqalign{\gamma(J_+)&\s d(\wedge\,\g)\cap( \wedge\,\g)^\g\cr &=0\cr}$$ But then (16) implies \vs {\bf Theorem 8.} {\it For any $k\in \Bbb Z$ one has $$R^k(\g)\s 
 H$$}\vs 

   Let $Sym(2k,2)$ be the subgroup of the symmetric group $Sym(2k)$ defined by $$\eqalign{&\,\,\,\,\,\,\,\,\,\,Sym(2k,2) = \cr &\{\sigma\in Sym(2k)\mid \sigma\,\,\hbox{permutes the set of unordered pairs}\,\,\{(1,2),(3,4),\ldots, ((2k-1),2k)\}\}\cr}$$
 That is if $\sigma\in Sym(2k,2)$ and $1\leq i\leq k$ there exists $1\leq j \leq k$ such that as unordered sets $$(\sigma(2i-1),\sigma (2i)) = ((2j-1),2j)$$\vskip.5pc It is clear that $Sym(2k,2)$ is a subgroup of order $2^k\cdot k!$. Let $\Pi(k)$ be a cross-section of the set of left cosets of $Sym(2k,2)$ in $Sym(2k)$ so that one has a disjoint union $$Sym(2k) = \cup\,\nu\,Sym(2k,2)\eqno (24)$$  indexed by   $\nu\in \Pi(k)$. \vs {\bf Remark 3}. One notes that the cardinality of $\Pi(k)$ is $(2k-1)(2k-3)\cdots 1$ and  the correspondence $$\nu \mapsto ((\nu(1),\nu(2)), (\nu(3),\nu(4)),\ldots ,(\nu((2k-1)),\nu(2k))) $$ sets up a bijection of $\Pi(k)$ with the set of all partitions of $(1,2,\ldots,2k)$ into a union of subsets each  of which has two elements. We also observe that $\Pi(k)$ may be chosen - and will be chosen -  such that $sg\,\nu =1$ for all $\nu\in \Pi(k)$. This is clear since the $sg$ character is not trivial on $Sym(k,2)$ for $k\geq 1$. The following is a restatement of the results in \S3.2 of [4] (see especially (3.25) and (3.29) in [4])\vs {\bf Theorem 9.} {\it For any $k\in \Bbb Z$ there exists a nonzero scalar $c_k$ such that for any  $x_i\,i=1,\dots, 2k,$ in $\g$
 $$\Gamma (x_1\wedge \cdots \wedge x_{2k}) = c_k \sum_{\nu\in \Pi(k)} [x_{\nu(1)},x_{\nu(2)}]\cdots [x_{\nu(2k-1)},x_{\nu(2k)}] \eqno (25))$$
 Furthermore the homogeneous polynomial of degree $k$ on the right side of (25) is harmonic and $R^k(\g)$ is the span of all such polynomials for an arbitrary choice of the $x_i$.} \vskip 1.pc  \centerline{\bf 4. On the variety of singular elements -- joint work with Nolan Wallach}\vskip 1.5pc Let $\hh$ be a Cartan sublgebra of $\g$ and let $\ell = dim\,\hh $ so $\ell =rank\,\g$. Let $\Delta$ be the set of roots of $(\hh,\g)$ and let $\Delta_+\s \Delta$ be a choice of 
 positive roots. Let $r = card\, \Delta_+$ so that $n = \ell + 2r$ where we fix $n= dim\,\g$. We assume a well ordering is defined on $\Delta_+$. For any $\varphi\in \Delta$ let $e_{\varphi}$ be a corresponding root vector. The choices will be normalized only insofar as $(e_{\varphi},e_{-\varphi}) = 1$ for all $\varphi\in \Delta$. From Proposition 2 one recovers  the well known fact that $\g^{(2k)}= 0$ for $k>r$ and $\g^{(2r)}$ is the set of all regular elements in $\g$. One also notes then that (16) implies $Var\, R^r(\g)$ reduces to $0$ if $k>r$ whereas Theorem 6 implies 
 $$Var\, R^{r}(\g)\,\,\hbox{is the set of all singular elements in $\g$}\eqno (26)$$ The paper [5]  is mainly devoted to a study of a special construction of $R^r(\g)$ and  a determination of its remarkable  $\g$-module structure.

It is a classic theorem of C. Chevalley that  $J$ is a polynomial ring in $\ell$ homogeneous generators $p_i$ so that we can write $$J = \Bbb C[p_1,\ldots,p_{\ell}]$$ Let $d_i = deg\,p_i$. Then if we put  $m_i = d_i-1$ the $m_i$ are referred to as the exponents of $\g$ and one knows $$\sum_{i=1}^{\ell} \,m_i = r\eqno (27)$$\vskip .5pc Henceforth assume $\g$ is simple so that the adjoint representation is irreducible. Let $y_j,\,j=1,\ldots,n,$ be  basis of $\g$. One defines an $\ell \times n$ matrix $Q = Q_{ij}, \,i=1,\ldots,\ell, \,j=1,\ldots,n$ by putting $$Q_{ij}=   \partial_{y_j}p_i\eqno (28)$$ Let $S_i,\, i=1,\ldots,\ell$, be the span of the entries of $Q$ in the $i^{th}$ row. The following is immediate \vs {\bf Proposition 3. } {\it $S_i\s P^{m_i}(\g)$. Furthermore $S_i$ is stable under the action of $\g$ and as a $\g$- module $S_i$ transforms according to the adjoint representation. }\vs If $V$ is a $\g$-module let $V_{ad}$ be the set of all of vectors in $V$ which transform according to the adjoint representation. The equality (24) readily implies $P(\g)_{ad} = J\otimes H_{ad}$.  I proved the following result some time ago (See \S5.4 in [2]. Especially see (5.4.6) and (5.4.7) in \S 5.4 of [2]) \vs  {\bf Theorem 10. } {\it The multiplicity of the adjoint representation in  $H_{ad}$ is $\ell$. Furthermore the invariants $p_i$ can be chosen so that $S_i\s H_{ad}$ for all $i$ and the $S_i, i=1,\ldots,\ell,$ are indeed the $\ell$ occurrences of the adjoint representation in $H_{ad}$.}\vs Clearly there are  $n\choose \ell$  $\ell\times \ell$ minors in the matrix $Q$. The determinant of any of these minors is an element of $P^r(\g)$ by (27). In [5] we offer a different formulation  of $R^r(\g)$ by proving \vs {\bf Theorem 11.} {\it The determinant of any $\ell\times\ell$ minor of $Q$ is an element of $R^r(\g)$ and indeed $R^r(\g) $ is the span of the determinants of all these minors. } \vskip 1.5pc \centerline{\bf 5. The $\g$ module structure of $R^ r(\g)$}\vskip 1.5pc The adjoint action of $\g$ on $\wedge\g $ extends to $U(\g)$ so that  $\wedge \g$ is a $U(\g)$-module. If $\ss\s \g $ is any subpace and $k= dim\,\ss$ let $[\ss]= 
 \wedge^k\ss$ so that $[\ss] $ is a 1-dimensional subspace of $\wedge^k\g$. Let $M_k\s \wedge^k\g $ be the span of all $[\ss]$ where $\ss$ is any $k$ dimensional commutative Lie subalgebra of $\g$. If no such subalgebra exists put $M_k=0$. It is clear that $M_k$ is a $\g$-submodule of $\wedge^k\g$. Let $Cas \in Z$ be the Casimir element corresponding to the Killing form. The following theorem was proved as Theorem (5) in [3] . \vs {\bf Theorem 12.} {\it For any $k\in \Bbb Z$ let $m_k$ be the maximal eigenvalue of $Cas$ on $\wedge^k\g$. Then $m_k\leq k$. Moreover $m_k=k$ if and only if $M_k\neq 0$ in which case $M_k$ is the eigenspace for the maximal eigenvalue $k$. }\vs Let $\Phi$ be a subset of  $\Delta$. Let $k =  card\,\Phi$ and  write, in increasing order, $$\Phi =\{\varphi_1,\ldots,\varphi_k\} \eqno (29)$$ Let $$e_{\Phi} = e_{\varphi_1}\wedge\cdots\wedge e_{\varphi_k}$$ so that $e_{\Phi}\in \wedge^k\g$ is an ($\hh$) weight vector with weight $$\langle \Phi\rangle = \sum_{i=1}^{k} \varphi_i$$ 
 
 Let $\n$ be the Lie algebra spanned by $e_{\varphi}$ for $\varphi\in \Delta_+$ and let $\b$ the Borel subalgebra of $\g$ defined by putting $\b= \hh + \n$. Now  a subset $\Phi\s \Delta_+$ will be called an ideal in $\Delta_+$ if the span, $\n_{\Phi}$, of $e_{\varphi}$, for $\varphi\in \Phi$, is an ideal of $\b$. In such a case $\Bbb C e_{\Phi}$ is stable under the action of $\b$ and hence if  $V_{\Phi} = U(\g)\cdot e_{\Phi}$ then, where $k = card \,\Phi$, $$ V_{\Phi}\s \wedge^k\g$$ is an irreducible $\g$-module of highest weight $\langle \Phi\rangle$ having $\Bbb C e_{\Phi}$ as the highest weight space.  We will say $\Phi$ is abelian if $\n_{\Phi}$ is an abelian ideal of $\b$. Let $${\cal A}(k)= \{\Phi\mid \hbox{$\Phi$ is an abelian ideal of cardinality $k$ in $\Delta_+$}\}$$\vskip .5pc The following theorem was established [3]. (See especially Theorems (7) and (8) in [3].) \vs {\bf Theorem 13.} {\it If $\Phi,\Psi$ are distinct ideals in $\Delta_+$ then $V_{\Phi}$
and $V_{\Psi}$ are inequivalent (i.e.$ \langle \Phi \rangle \neq  \langle \Psi \rangle$). Furthermore if $M_k\neq 0$ then $$M_k = \oplus_{\Phi\in {\cal A}(k)} V_{\Phi}\eqno (30)$$ so that, in particular, $M_k$ is a multiplicity 1 $\g$-module. }\vs We now focus on the case where $k=\ell$. Clearly $M_{\ell}\neq 0$ since $\g^x$ is an abelian subalgebra of dimension $\ell$ for any regular $x\in \g$. Let ${\cal I}(\ell)$ be the set of all ideals of cardinality $\ell$. The following theorem giving the remarkable structure of $R^r(\g)$ as a $\g$-module is one of the main results in [4]]. \vs {\bf Theorem 14.} {\it One has ${\cal I}(\ell) = {\cal A}(\ell)$ so that $$M_{\ell} =  \oplus _{\Phi\in {\cal I}(\ell)} V_{\Phi}\eqno (31)$$ Moreover as 
$\g$-modules one has the equivalence $$R^r(\g)\cong M_{\ell}\eqno (32)$$ so that $R^r(\g)$ is a multiplicity 1 $\g$-module with $card\,{\cal I}(\ell)$  irreducible components and $Cas$ takes the value $\ell$ on each and every one of the ${\cal I}(\ell)$ distinct components.} \vskip 1.5pc 
{\bf Example}. If $\g$ is of type $A_{\ell}$ then then the elements of ${\cal I}(\ell)$ can identified with Young diagrams of size $\ell$. In this case therefore the number of irreducible components in $R^{r}(\g)$ is $P(\ell)$ where $P$ here is the classical partition function. \vskip 1.5pc

 {\bf References}

\vskip 1.5pc

[1] B. Kostant, A Theorem of Frobenius, a Theorem of Amitsur-Levitski and Cohomology Theory, {\it  J. Mech and Mech.,} {\bf 7} (1958): 2, Indiana University, 237--264.

\vskip 1.5pc

[2]'' , Lie Group Representations on Polynomial Rings, {\it  American J. Math},

{\bf 85} (1963), No. 1, 327--404.

\vskip 1.pc

[3]'', Eigenvalues of a Laplacian and Commutative Lie Subalgebras,  {\it Topology}, {\bf 13}, (1965), 147--159.

\vskip 1.5pc

[4]'', A Lie Algebra Generalization of the Amitsur-Levitski Theorem,   {{ {\it Adv. In Math.}, {\bf 40} (1981):2, 155--175.

\vskip 1.5pc

[5]'' and N. Wallach,  On the algebraic set of singular elements in a complex simple Lie algebra, in: {\it Representation Theory and Mathematical Physics}, Conf. in honor of Gregg Zuckerman's 60th Birthday, Contemp. Math.,
{\bf  557},  Amer. Math. Soc., 2009,  pp. 215--230.

 \end

 \end

 \end

 \end

 \end

 \end